\documentclass[a4paper,11pt]{amsart}
\usepackage{graphicx}
\usepackage{mathptmx}
\usepackage{amsmath}
\usepackage{amssymb}
\usepackage{enumitem}
\usepackage{xcolor}

\newmuskip\pFqmuskip

\newcommand*\pFq[6][8]{%
  \begingroup % only local assignments
  \pFqmuskip=#1mu\relax
  \mathcode`=\string"8000
  \begingroup\lccode`\~=`\,
  \lowercase{\endgroup\let~}\pFqcomma
  F^{#2}_{#3}{\left(\genfrac..{0pt}{}{#4}{#5}\bigg|#6\right)}%
  \endgroup
}
\newcommand{\pFqcomma}{\mskip\pFqmuskip}

\newtheorem{theorem}{Theorem}
\newtheorem{lemma}[theorem]{Lemma}

\begin{document}

\title[]{Two variable degenerate Bell polynomials associated with Poisson degenerate central moments}

\author{Dae San Kim}
\address{Department of Mathematics, Sogang University, Seoul 121-742, Republic of Korea}
\email{dskim@sogang.ac.kr}

\author{Taekyun  Kim}
\address{Department of Mathematics, Kwangwoon University, Seoul 139-701, Republic of Korea}
\email{tkkim@kw.ac.kr}

\author{Hanyoung Kim}
\address{Department of Mathematics, Kwangwoon University, Seoul 139-701, Republic of Korea}
\email{gksdud213@kw.ac.kr}

\author{Hyunseok Lee}
\address{Department of Mathematics, University, Seoul 139-701, Republic of Korea}
\email{luciasconstant@kw.ac.kr}

\subjclass[2010]{11B73; 11B83}
\keywords{two variable degenerate Bell polynomials; Poisson degenerate central moments; Stirling type numbers}

\maketitle

\begin{abstract}
In this paper, we introduce two variable degenerate Bell polynomials as a generalization of the degenerate Bell polynomials which were studied earlier. Among other things, we show their connections with Poisson degenerate central moments and Charlier polynomials.
\end{abstract}

\pagestyle{myheadings}
\markboth{\centerline{\scriptsize D.S. Kim, T. Kim, H. Kim, H. Lee}}
          {\centerline{\scriptsize Two variable degenerate Bell polynomials associated with Poisson degenerate central moments}}
 
\section{Introduction} 

In recent years, the idea of studying degenerate versions of many special polynomials and numbers, which was initiated by Carlitz in [4], received regained interests of some mathematicians and many interesting results were discovered (see\,\,[7,8,10,11,13-15]). They have been explored by employing several different tools such as combinatorial methods, generating functions, $p$-adic analysis, umbral calculus techniques, differential equations and probability theory. \par
The aim of this paper is to introduce two variable degenerate Bell polynomials  and study their connections with Poisson degenerate central moments and Charlier polynomials. Here we note that those polynomials are  generalizations of the single variable degenerate Bell polynomials and in turn these single variable degenerate Bell polynomials are degenerate versions of the single variable Bell polynomials (see\,\,[11,13,15]). \par
The outline of this paper is as follows. In Section 1, as a preparation to the next section we will recall the degenerate exponential functions, the Stirling numbers of the first and second kinds, Bell polynomials and numbers, the degenerate Stirling numbers of the first and second kinds, the Poisson random variable with parameter $\alpha$, the degenerate Bell polynomials and numbers, and the Charlier polynomials. In Section 2, we will introduce the two variable degenerate Bell polynomials, express them in terms of the degenerate Bell polynomials and show they reduce to the degenerate Bell polynomials when the two arguments are the same. We will define the Poisson degenerate central moments and prove that they are equal to the two variable degenerate Bell polynomials. Also, we will find an explicit expression for the two variable degenerate Bell polynomials. By letting one of the arguments of the two variable degenerate Bell polynomials equal to zero, we will be led to define Stirling type numbers as the coefficients in the expansion of powers of $\alpha$. Then we will express the Stirling type numbers in terms of the degenerate Stirling numbers of the second kind and vice versa. Finally, we will be able to find certain connections between the two variable degenerate Bell polynomials and Charlier polynomials involving degenerate Stirling numbers. \par

\vspace{0.1in}

For any $0 \neq\lambda\in\mathbb{R}$, the degenerate exponential functions are defined by 
\begin{equation}
	e_{\lambda}^{x}(t)=(1+\lambda t)^{\frac{x}{\lambda}}=\sum_{n=0}^{\infty}(x)_{n,\lambda}\frac{t^{n}}{n!},\quad(\mathrm{see}\,\, [7,13,14]), \label{1}
\end{equation}
where $(x)_{0,\lambda}=1$, $(x)_{n,\lambda}=x(x-\lambda)(x-2\lambda)\cdots\big(x-(n-1)\lambda\big)$, $(n\ge 1)$. \\
In \eqref{1}, we note that $\displaystyle \lim_{\lambda\rightarrow 0}e_{\lambda}^{x}(t)=e^{tx}\displaystyle$, and let $e_{\lambda}(t)=e_{\lambda}^{1}(t)$. \par
As is well known, the Stirling numbers of the first kind are defined as 
\begin{equation}
	(x)_{n}=\sum_{l=0}^{n}S_{1}(n,l)x^{l},\quad(n\ge 0),\quad(\mathrm{see}\,\, [8,10]), \label{2}
\end{equation}
where $(x)_{0}=1$, $(x)_{n}=x(x-1)(x-2)\cdots\big(x-(n-1)\big)$, $(n\ge 1)$. 
As the inversion formula of \eqref{2}, the Stirling numbers of the second kind are defined by 
\begin{equation}
	x^{n}=\sum_{l=0}^{n}S_{2}(n,l)(x)_{l},\quad(n\ge 0),\quad(\mathrm{see}\,\, [8,10]). \label{3}
\end{equation}
From \eqref{2} and \eqref{3}, we note that 
\begin{equation}
	\frac{1}{k!}\big(\log(1+t)\big)^{k}=\sum_{n=k}^{\infty}S_{1}(n,k)\frac{t^{n}}{n!}, \label{4}
\end{equation}
and 
\begin{equation}
	\frac{1}{k!}\big(e^{t}-1\big)^{k}=\sum_{n=k}^{\infty}S_{2}(n,k)\frac{t^{n}}{n!},\quad(k\ge 0).\label{5}
\end{equation} \par
The Bell polynomials are given by  
\begin{equation}
	e^{x(e^{t}-1)}=\sum_{n=0}^{\infty}\mathrm{Bel}_{n}(x)\frac{t^{n}}{n!},\quad(\mathrm{see}\,\, [1-3,5,6,17]).\label{6}
\end{equation}
Thus, we have 
\begin{equation}
\mathrm{Bel}_{n}(x)=\sum_{l=0}^{n}S_{2}(n,l)x^{l}=e^{-x}\sum_{k=0}^{\infty}\frac{k^n}{k!}x^k,\quad(\mathrm{see}\,\, [1-3,5,6,17]). \label{7}
\end{equation}
When $x=1$, $\mathrm{Bel}_{n}= \mathrm{Bel}_{n}(1)$ are called the Bell numbers. \par
The degenerate Stirling numbers of the second kind are defined by 
\begin{equation}
	(x)_{n,\lambda}=\sum_{l=0}^{n}S_{2,\lambda}(n,l)(x)_{l},\quad(n\ge 0),\quad(\mathrm{see}\,\, [7,10,13-15]). \label{8}
\end{equation}
As the inversion formula of \eqref{8}, the degenerate Stirling numbers of the first kind are given by 
\begin{equation}
	(x)_{n}=\sum_{l=0}^{n}S_{1,\lambda}(n,l)(x)_{l,\lambda},\quad(n\ge 0),\quad(\mathrm{see}\,\, [7,14]). \label{9}
\end{equation}
Let $\log_{\lambda}(t)$ be the compositional inverse of $e_{\lambda}(t)$ such that $e_{\lambda}(\log_{\lambda}t)=\log_{\lambda}\big(e_{\lambda}(t)\big)=t$.  
Then we note that 
\begin{equation}
\log_{\lambda}(1+t)=\frac{1}{\lambda}\big((1+t)^{\lambda}-1\big)=\sum_{n=1}^{\infty}\lambda^{n-1}(1)_{n,\frac{1}{\lambda}}\frac{t^{n}}{n!},\quad(\mathrm{see}\,\, [7,14]). \label{10}
\end{equation}
From \eqref{10}, we note that $\displaystyle\lim_{\lambda\rightarrow 0}\log_{\lambda}(1+t)=\log(1+t) \displaystyle$. 
By \eqref{8} and \eqref{9}, we easily get 
\begin{equation}
\frac{1}{k!}\big(e_{\lambda}(t)-1\big)^{k}=\sum_{n=k}^{\infty}S_{2,\lambda}(n,k)\frac{t^{n}}{n!},\quad(\mathrm{see}\,\, [7,10,13-15]), \label{11}
\end{equation}
and 
\begin{equation}
\frac{1}{k!}\big(\log_{\lambda}(1+t)\big)^{k}=\sum_{n=k}^{\infty}S_{1,\lambda}(n,k)\frac{t^{n}}{n!},\quad(\mathrm{see}\,\, [7,14]), \label{12}
\end{equation}
where $k$ is a nonnegative integer. \par 
A random variable $X$ is a real valued function defined on a sample space. If $X$ takes any values in a countable set, then $X$ is called a discrete random variable. For a discrete random variable $X$, the probability mass function $p(a)$ of $X$ is defined by 
\begin{displaymath}
	p(a)=P\{X=a\},\quad(\mathrm{see}\,\, [13,19]). 
\end{displaymath}
A random variable $X$ taking on one of the values $0,1,2,\dots$ is said to be the Poisson random variable with parameter $\alpha(>0)$ if the probability mass function of $X$ is given by 
\begin{equation}
	p(i)=P\{X=i\}=e^{-\alpha}\frac{\alpha^{i}}{i!},\quad i=0,1,2,\dots. \label{13}
\end{equation}
For $n\ge 1$, the quantity $E[X^{n}]$ of the Poisson random variable $X$ with parameter $\alpha(>0)$, which is called the $n$-th moment of $X$, is given by 
\begin{equation}
	E[X^{n}]=\sum_{i=0}^{\infty}i^{n}p(i)=e^{-\alpha}\sum_{i=0}^{\infty}\frac{i^{n}}{i!}\alpha^{i}=\mathrm{Bel}_{n}(\alpha). \label{14} 
\end{equation} \par
Recently, Kim-Kim-Dolgy introduced the degenerate Bell polynomials given by
\begin{equation}
	e^{x(e_{\lambda}(t)-1)}=\sum_{n=0}^{\infty}\mathrm{Bel}_{n,\lambda}(x)\frac{t^{n}}{n!},\quad(\mathrm{see}\,\, [11,13,15]).\label{15}
\end{equation}
When $x=1$, $\mathrm{Bel}_{n,\lambda}= \mathrm{Bel}_{n,\lambda} (1)$ are called the degenerate Bell numbers. 
It is not difficult to show that $\displaystyle \lim_{\lambda\rightarrow 0}\mathrm{Bel}_{n,\lambda}(x)=\mathrm{Bel}_{n}(x) \displaystyle$.\\ 
From \eqref{15}, we note that 
\begin{align}
\mathrm{Bel}_{n,\lambda}(x)\ =\ \sum_{l=0}^{n}x^{l}S_{2,\lambda}(n,l) \label{16} 
=\ 	e^{-x}\sum_{k=0}^{\infty}\frac{(k)_{n,\lambda}}{k!}x^{k},\quad(\mathrm{see}\,\, [11,13,15]). \end{align} \par
It is well known that Charlier polynomials $C_{n}(x;\alpha)$ are defined by 
\begin{equation}
	e^{-\alpha t}(1+t)^{x}=\sum_{n=0}^{\infty}C_{n}(x;\alpha)\frac{t^{n}}{n!},\quad(\mathrm{see}\,\, [16-18]), \label{17}
\end{equation}
where $x,t,\alpha\in\mathbb{R}$. 
We observe that
\begin{align}
e^{-\alpha t}(1+t)^{x}\ & =\ \sum_{m=0}^{\infty}(-\alpha)^{m}\frac{t^{m}}{m!}\sum_{l=0}^{\infty}x^{l}\frac{\big(\log(1+t)\big)^{l}}{l!} \label{18} \\
&=\ \sum_{m=0}^{\infty}(-\alpha)^{m}\frac{t^{m}}{m!}\sum_{k=0}^{\infty}\sum_{l=0}^{k}x^{l}S_{1}(k,l)\frac{t^{k}}{k!} \nonumber \\
&=\ \sum_{n=0}^{\infty}\bigg(\sum_{k=0}^{n}(-\alpha)^{n-k}\binom{n}{k}\sum_{l=0}^{k}x^{l}S_{1}(k,l)\bigg)\frac{t^{n}}{n!}. \nonumber
\end{align}
Thus, by \eqref{17} and \eqref{18}, we get 
\begin{align}\label{19}
C_{n}(x;\alpha)\  
	=\ \sum_{l=0}^{n}\bigg(\sum_{k=l}^{n}\binom{n}{k}(-1)^{n-k}\alpha^{n-k}S_{1}(n,k)\bigg)x^{l},\quad (n\ge 0). 
\end{align} \par

\section{Two Variable Degenerate Bell polynomials} 
For $\alpha,\beta\in\mathbb{R}$, let us consider {\it{two variable degenerate Bell polynomials}} which are given by
\begin{equation}
e^{\alpha(e_{\lambda}(t)-1)}e_{\lambda}^{\beta-\alpha}(t)=\sum_{n=0}^{\infty}\mathrm{Bel}_{n,\lambda}(\alpha,\beta)\frac{t^{n}}{n!}. \label{20}
\end{equation}
From the left hand side of \eqref{20}, we have 
\begin{align}
e^{\alpha(e_{\lambda}(t)-1)}e_{\lambda}^{\beta-\alpha}(t)\ &=\ \sum_{l=0}^{\infty}\mathrm{Bel}_{l,\lambda}(\alpha)\frac{t^{l}}{l!}\sum_{m=0}^{\infty}(\beta-\alpha)_{m,\lambda}\frac{t^{m}}{m!} \nonumber \\
	&=\ \sum_{n=0}^{\infty}\bigg(\sum_{l=0}^{n}\binom{n}{l}\mathrm{Bel}_{l,\lambda}(\alpha)(\beta-\alpha)_{n-l,\lambda}\bigg)\frac{t^{n}}{n!}. \label{21}
\end{align}
Therefore, by \eqref{20} and \eqref{21}, we obtain the following theorem. 
\begin{theorem}
	For $n\ge 0$, we have 
	\begin{displaymath}
	\mathrm{Bel}_{n,\lambda}(\alpha,\beta)=\sum_{l=0}^{n}\binom{n}{l}\mathrm{Bel}_{l,\lambda}(\alpha)(\beta-\alpha)_{n-l,\lambda}
	\end{displaymath}
	In particular, 
	\begin{displaymath}
	\mathrm{Bel}_{n,\lambda}(\alpha,\alpha)= \mathrm{Bel}_{n,\lambda}(\alpha),\quad(n\ge 0). 
	\end{displaymath}
\end{theorem}
Let $X$ be a Poisson random variable with parameter $\alpha(>0)$. 
Then we consider the {\it{Poisson degenerate central moments}} given by $E[(X-\alpha)_{n,\lambda}]$, $(n\ge 0)$, \\
where $(X-\alpha)_{0,\lambda}=1$, $(X-\alpha)_{n,\lambda}=(X-\alpha)(X-\alpha-\lambda)\cdots(X-\alpha-(n-1)\lambda)$, $(n\ge 1)$. \\
Note that 
\begin{displaymath}
	\lim_{\lambda\rightarrow 0}E\big[(X-\alpha)_{n,\lambda}\big]=E\big[(X-\alpha)^{n}\big].
\end{displaymath} \par
We observe that 
\begin{align}
E\big[e_{\lambda}^{X-\alpha+\beta}(t)\big]\ &=\ \sum_{k=0}^{\infty}e_{\lambda}^{k-\alpha+\beta}(t)p(k)\label{22}\\
&=\ e_{\lambda}^{\beta-\alpha}(t)e^{-\alpha}\sum_{k=0}^{\infty}e_{\lambda}^{k}(t) \frac{\alpha^{k}}{k!}\nonumber \\
&=\ e_{\lambda}^{\beta-\alpha}(t)e^{\alpha(e_{\lambda}(t)-1)}\nonumber \\ 
&=\  \sum_{n=0}^{\infty}\mathrm{Bel}_{n,\lambda}(\alpha,\beta)\frac{t^{n}}{n!}. \nonumber	
\end{align}
On the other hand, 
\begin{align}
	E\big[e_{\lambda}^{X-\alpha+\beta}(t)\big]\ &=\ E\bigg[\sum_{n=0}^{\infty}(X-\alpha+\beta)_{n,\lambda}\frac{t^{n}}{n!}\bigg] \label{23} \\
	&=\ \sum_{n=0}^{\infty}E\big[(X-\alpha+\beta)_{n,\lambda}\big]\frac{t^{n}}{n!}. \nonumber
\end{align}
Therefore, by \eqref{22} and \eqref{23}, we obtain the following theorem. 
\begin{theorem}
	For $n\ge 0$, we have 
	\begin{displaymath}
	E\big[(X-\alpha+\beta)_{n,\lambda}\big]\ =\ \mathrm{Bel}_{n,\lambda}(\alpha,\beta).
	\end{displaymath}
	In particular, 
	\begin{displaymath}
	E\big[(X-\alpha+\beta)^{n}\big]\ =\ \sum_{l=0}^{n}\binom{n}{l}\mathrm{Bel}_{l}(\alpha)(\beta-\alpha)^{n-l}. 
	\end{displaymath}
\end{theorem}
It is easy to show that 
\begin{equation}
	(x+y)_{n,\lambda}\ =\ \sum_{i=0}^{n}\binom{n}{i}(x)_{i,\lambda}(y)_{n-i,\lambda},\quad (n\ge 0). \label{24}
\end{equation}
By \eqref{24}, we get 
\begin{equation}
	E\big[(X-\alpha+\beta)_{n,\lambda}\big]\ =\ \sum_{k=0}^{n}\binom{n}{k}(\beta-\alpha)_{n-k,\lambda}E\big[(X)_{k,\lambda}\big].\label{25}
	\end{equation} \par
From \eqref{13} and \eqref{15}, we have 
\begin{align}
E\big[e_{\lambda}^{X}(t)\big]\ &=\ \sum_{n=0}^{\infty}e_{\lambda}^{n}(t)p(n)\ =\ e^{-\alpha}\sum_{n=0}^{\infty}e_{\lambda}^{n}(t)\frac{\alpha^{n}}{n!}\label{26} \\
&=\ e^{\alpha(e_{\lambda}(t)-1)}\ =\ \sum_{n=0}^{\infty}\mathrm{Bel}_{n,\lambda}(\alpha)\frac{t^{n}}{n!}.\nonumber
\end{align}
On the other hand, 
\begin{equation}
E\big[e_{\lambda}^{X}(t)\big]=E\bigg[\sum_{n=0}^{\infty}(X)_{n,\lambda}\frac{t^{n}}{n!}\bigg]=\sum_{n=0}^{\infty}E[(X)_{n,\lambda}]\frac{t^{n}}{n!}. \label{27}
\end{equation}
By \eqref{26} and \eqref{27}, we get 
\begin{equation}
	\mathrm{Bel}_{n,\lambda}(\alpha)\ =\ E[(X)_{n,\lambda}],\quad(n\ge 0). \label{28}
\end{equation} \par
From Theorem 2, \eqref{25} and \eqref{28}, we have 
\begin{align}
\mathrm{Bel}_{n,\lambda}(\alpha,\beta)\ &=\ E\big[(X-\alpha+\beta)_{n,\lambda}\big]\ =\ \sum_{k=0}^{n}\binom{n}{k}(\beta-\alpha)_{n-k,\lambda}E\big[(X)_{k,\lambda}\big] \nonumber \\
&=\ \sum_{k=0}^{n}\binom{n}{k}(\beta-\alpha)_{n-k,\lambda}\mathrm{Bel}_{k,\lambda}(\alpha)\nonumber\\
&=\ \sum_{k=0}^{n}\binom{n}{k}(\beta-\alpha)_{n-k,\lambda}\sum_{i=0}^{k}\alpha^{i}S_{2,\lambda}(k,i) \label{29} \\
&=\ \sum_{i=0}^{n}\alpha^{i}\sum_{k=i}^{n}\binom{n}{k}(\beta-\alpha)_{n-k,\lambda}S_{2,\lambda}(k,i),\quad(n\ge 0).\nonumber 	
\end{align}
Therefore, by \eqref{29}, we obtain the following theorem. 
\begin{theorem}
	For $n\ge 0$, we have 
	\begin{displaymath}
		\mathrm{Bel}_{n,\lambda}(\alpha,\beta)\ =\ \sum_{i=0}^{n}\alpha^{i}\sum_{k=i}^{n}\binom{n}{k}(\beta-\alpha)_{n-k,\lambda}S_{2,\lambda}(k,i).
	\end{displaymath}
\end{theorem}
Let us take $\beta=0$ in Theorem 2. Then we have 
\begin{equation}
	E\big[(X-\alpha)_{n,\lambda}]\ =\ \mathrm{Bel}_{n,\lambda}(\alpha,0),\quad(n\ge 0).\label{30}
\end{equation}
In light of \eqref{16}, we define {\it{Stirling type numbers}} $S_{\lambda}(n,k)$ associated with $\mathrm{Bel}_{n,\lambda}(\alpha,0)$ as
\begin{equation}
	E\big[(X-\alpha)_{n,\lambda}\big]\ =\ \mathrm{Bel}_{n,\lambda}(\alpha,0)\ =\ \sum_{k=0}^{n}S_{\lambda}(n,k)\alpha^{k},\quad (n\ge 0). \label{31}
\end{equation}
When $\alpha=1$, the Poisson degenerate central moments are given by 
\begin{equation}
	\mathrm{Bel}_{n,\lambda}(1,0)\ =\ E\big[(X-1)_{n,\lambda}\big]\ =\ \sum_{k=0}^{n}S_{\lambda}(n,k), \quad(n\ge 0). \label{32}
\end{equation}
By Theorem 2 and \eqref{24}, we easily get  
\begin{align}
\mathrm{Bel}_{n,\lambda}(\alpha,\beta)\ &=\ E\big[(X-\alpha+\beta)_{n,\lambda}\big]\ =\ \sum_{k=0}^{n}\binom{n}{k}(\beta)_{n-k,\lambda}E\big[(X-\alpha)_{k,\lambda}\big] \label{33} \\
&=\ \sum_{k=0}^{n}\binom{n}{k}(\beta)_{n-k,\lambda}\mathrm{Bel}_{k,\lambda}(\alpha,0)\ =\ \sum_{k=0}^{n}\binom{n}{k}(\beta)_{n-k,\lambda}\sum_{l=0}^{k}S_{\lambda}(k,l)\alpha^{l}\nonumber \\	
&=\ \sum_{l=0}^{n}\alpha^{l}\sum_{k=l}^{n}\binom{n}{k}(\beta)_{n-k}S_{k}(k,l). \nonumber
\end{align}
From \eqref{29}, we note that 
\begin{align}
\mathrm{Bel}_{n,\lambda}(\alpha,\beta)\ &=\ \sum_{k=0}^{n}\binom{n}{k}(\beta-\alpha)_{k,\lambda}\mathrm{Bel}_{n-k,\lambda}(\alpha)\label{34}	\\
&=\ \sum_{k=0}^{n}\binom{n}{k}\sum_{l=0}^{k}\binom{k}{l}(\beta)_{l,\lambda}(-\alpha)_{k-l,\lambda}\sum_{i=0}^{n-k}\alpha^{i}S_{2,\lambda}(n-k,i)\nonumber \\
&=\ \sum_{k=0}^{n}\sum_{l=0}^{k}\binom{n}{l}\binom{n-l}{n-k}(\beta)_{l,\lambda}(-\alpha)_{k-l,\lambda}\sum_{i=0}^{n-k}\alpha^{i}S_{2,\lambda}(n-k,i) \nonumber \\
&=\ \sum_{l=0}^{n}\sum_{b=0}^{n-l}\binom{n}{l}\binom{n-l}{b}(\beta)_{l,\lambda}(-\alpha)_{n-b-l,\lambda}\sum_{i=0}^{b}\alpha^{i}S_{2,\lambda}(b,i) \nonumber\\
&=\ \sum_{l=0}^{n}\sum_{b=0}^{l}\binom{n}{n-l}\binom{l}{b}(\beta)_{n-l,\lambda}(-\alpha)_{l-b,\lambda}\sum_{i=0}^{b}\alpha^{i}S_{2,\lambda}(b,i) \nonumber\\
&=\ \sum_{l=0}^{n}\sum_{b=0}^{l}\binom{n}{l}\binom{l}{l-b}(\beta)_{n-l,\lambda}(-\alpha)_{b,\lambda}\sum_{i=0}^{l-b}\alpha^{i}S_{2,\lambda}(l-b,i) \nonumber \\
&=\ \sum_{l=0}^{n}\sum_{b=0}^{l}\binom{n}{l}\binom{l}{b}(\beta)_{n-l,\lambda}(-\alpha)_{b,\lambda}\sum_{i=b}^{l}\alpha^{i-b}S_{2,\lambda}(l-b,i-b) \nonumber\\
&=\ \sum_{l=0}^{n}\binom{n}{l}(\beta)_{n-l,\lambda}\sum_{i=0}^{l}\alpha^{i}\sum_{b=0}^{i}\binom{l}{b}(-1)^{b}\langle 1\rangle_{b,\frac{\lambda}{\alpha}}S_{2,\lambda}(l-b,i-b), \nonumber
\end{align}
where $\langle x\rangle_{0,\lambda}=1$, $\langle x\rangle_{n,\lambda}=x(x+\lambda)(x+2\lambda)\cdots\big(x+(n-1)\lambda\big)$, $(n\ge 1)$. \\
Thus, by \eqref{34}, we get 
\begin{align}
	\mathrm{Bel}_{n,\lambda}(\alpha,0)\ &=\ \sum_{i=0}^{n}\alpha^{i}\sum_{b=0}^{i}\binom{n}{b}(-1)^{b}\langle 1\rangle_{b,\frac{\lambda}{\alpha}}S_{2,\lambda}(n-b,i-b)\label{35}\\
	&=\ E\big[(X-\alpha)_{n,\lambda}\big],\quad(n\ge 0). \nonumber
\end{align} 
Therefore, by \eqref{31} and \eqref{35}, we obtain the following theorem. 
\begin{theorem}
	For $n,i\ge 0$, we have 
	\begin{displaymath}
	S_{\lambda}(n,i)\ =\ \sum_{b=0}^{i}\binom{n}{b}(-1)^{b}\langle 1\rangle_{b,\frac{\lambda}{\alpha}}S_{2,\lambda}(n-b,i-b), 
	\end{displaymath}
where $\langle x\rangle_{0,\lambda}=1$, $\langle x\rangle_{x,\lambda}=x(x+\lambda)\cdots\big(x+(n-1)\lambda\big)$, $(n\ge 1)$. 
\end{theorem}
Now, we consider the inversion formula of Theorem 4. We observe that 
\begin{equation}
	\sum_{k=0}^{n}\alpha^{k}S_{2,\lambda}(n,k)\ =\ \mathrm{Bel}_{n,\lambda}(\alpha)\ =\ \mathrm{Bel}_{n,\lambda}(\alpha,\alpha)\ =\ E\big[(X-\alpha+\alpha)_{n,\lambda}\big]\label{36}
\end{equation}
\begin{align*}
	&=\ \sum_{l=0}^{n}\binom{n}{l}(\alpha)_{n-l,\lambda}E\big[(X-\alpha)_{l,\lambda}\big]\ =\ \sum_{l=0}^{n}\binom{n}{l}(\alpha)_{n-l,\lambda}\sum_{k=0}^{l}\alpha^{k}S_{\lambda}(l,k)  \\
	&=\ \sum_{l=0}^{n}\binom{n}{n-l}(\alpha)_{l,\lambda}\sum_{k=0}^{n-l}\alpha^{k}S_{\lambda}(n-l,k)\ =\ \sum_{l=0}^{n}\binom{n}{l}(\alpha)_{l,\lambda}\sum_{k=l}^{n}\alpha^{k-l}S_{\lambda}(n-l,k-l)\\
	&=\ \sum_{k=0}^{n}\alpha^{k}\sum_{l=0}^{k}\binom{n}{l}\alpha^{-l}(\alpha)_{l,\lambda}S_{\lambda}(n-k,k-l)\ =\ \sum_{k=0}^{n}\alpha^{k}\sum_{l=0}^{k}\binom{n}{l}(1)_{l,\frac{\lambda}{\alpha}}S_{\lambda}(n-l,k-l).
\end{align*}
Therefore, by comparing the coefficients on both sides of \eqref{36}, we obtain the following theorem. 
\begin{theorem}
	For $n,k\ge 0$, we have 
	\begin{displaymath}
		S_{2,\lambda}(n,k)\ =\ \sum_{l=0}^{k}\binom{n}{l}(1)_{l,\frac{\lambda}{\alpha}}S_{\lambda}(n-l,k-l).
	\end{displaymath}
\end{theorem}
We note that the left hand side of \eqref{17} is given by
\begin{equation}
	e^{-\alpha t}(1+t)^{x}\ =\ e^{-\alpha t}e_{\lambda}^{x}\big(\log_{\lambda}(1+t)\big)\ =\ e^{-\alpha t}\sum_{l=0}^{\infty}(x)_{l,\lambda}\frac{1}{l!}\big(\log_{\lambda}(1+t)\big)^{l} \label{37}
\end{equation}
\begin{align*}
&=\	e^{-\alpha t}\sum_{l=0}^{\infty}(x)_{l,\lambda}\sum_{k=l}^{\infty}S_{1,\lambda}(k,l)\frac{t^{k}}{k!}\ =\ \sum_{m=0}^{\infty}(-\alpha)^{m}\frac{t^m}{m!}\sum_{k=0}^{\infty}\sum_{l=0}^{k}(x)_{l,\lambda}S_{1,\lambda}(k,l)\frac{t^{k}}{k!} \\
&=\ \sum_{n=0}^{\infty}\bigg(\sum_{k=0}^{n}\binom{n}{k}(-\alpha)^{n-k}\sum_{l=0}^{k}(x)_{l,\lambda}S_{1,\lambda}(k,l)\bigg)\frac{t^{n}}{n!}. 
\end{align*}
Therefore, by \eqref{17} and \eqref{37}, we obtain the following lemma. 
\begin{lemma}
	For $n\ge 0$, we have 
	\begin{displaymath}
		C_{n}(x;\alpha)\ =\ \sum_{l=0}^{n}(x)_{l,\lambda}\sum_{k=l}^{n}\binom{n}{k}(-\alpha)^{n-k}S_{1,\lambda}(k,l).
	\end{displaymath}
\end{lemma}
Replacing $t$ by $e_{\lambda}(t)-1$ and $x$ by $-\alpha$ in \eqref{17}, we get 
\begin{align}
e^{\alpha(e_{\lambda}(t)-1)}e_{\lambda}^{-\alpha}(t)\ &=\ \sum_{k=0}^{\infty}C_{k}\big(-\alpha;-\alpha\big)\frac{1}{k!}\big(e_{\lambda}(t)-1\big)^{k}\label{38} \\
&=\ \sum_{k=0}^{\infty}C_{k}(-\alpha;-\alpha)\sum_{n=k}^{\infty}S_{2,\lambda}(n,k)\frac{t^{n}}{n!}\nonumber \\
&=\ \sum_{n=0}^{\infty}\bigg(\sum_{k=0}^{n}C_{k}(-\alpha,-\alpha)S_{2,\lambda}(n,k)\bigg)\frac{t^{n}}{n!}.\nonumber
\end{align}
On the other hand, 
\begin{equation}
e^{\alpha(e_{\lambda}(t)-1)}e_{\lambda}^{-\alpha}(t)\ =\ \sum_{n=0}^{\infty}\mathrm{Bel}_{n,\lambda} (\alpha,0)\frac{t^{n}}{n!}. \label{39}
\end{equation}
In particular, for $x=0$ in \eqref{17}, we get 
\begin{displaymath}
e^{-\alpha t}\ =\ \sum_{k=0}^{\infty}C_{k}(0;\alpha)\frac{t^{k}}{k!}.
\end{displaymath}
By replacing $t$ by $1-e_{\lambda}(t)$, we get 
\begin{align}
e^{\alpha(e_{\lambda}(t)-1)}\ &= \sum_{k=0}^{\infty}C_{k}(0;\alpha)\frac{1}{k!}\big(1-e_{\lambda}(t)\big)^{k} \label{40} \\
&=\ \sum_{k=0}^{\infty}(-1)^{k}C_{k}(0;\alpha)\sum_{n=k}^{\infty}S_{2,\lambda}(n,k)\frac{t^{n}}{n!} \nonumber \\
&=\ \sum_{n=0}^{\infty}\bigg(\sum_{k=0}^{n}(-1)^{k}C_{k}(0;\alpha)S_{2,\lambda}(n,k)\bigg)\frac{t^{n}}{n!}\nonumber
\end{align}
Therefore, by \eqref{38}, \eqref{39} and \eqref{40}, we obtain the following lemma. 
\begin{lemma}
	For $n\ge 0$, we have 
	\begin{displaymath}
	\mathrm{Bel}_{n,\lambda}(\alpha,0)=\sum_{k=0}^{n}C_{k}(-\alpha;-\alpha)S_{2,\lambda}(n,k),
	\end{displaymath}
	and
	\begin{displaymath}
	\mathrm{Bel}_{n,\lambda}(\alpha)=\sum_{k=0}^{n}(-1)^{k}C_{k}(0;\alpha)S_{2,\lambda}(n,k).
	\end{displaymath}
\end{lemma}
From \eqref{17}, we get 
\begin{equation}
	e^{\alpha t}(1+t)^{\beta-\alpha}\ =\ \sum_{n=0}^{\infty}C_{n}(\beta-\alpha\ ;\ -\alpha)\frac{t^{n}}{n!}.\label{41}
\end{equation}
Replacing $t$ by $e_{\lambda}(t)-1$ in \eqref{41}, we get 
\begin{align}
	e^{\alpha(e_{\lambda}(t)-1)}e_{\lambda}^{\beta-\alpha}(t)\ &=\ \sum_{k=0}^{\infty}C_{k}(\beta-\alpha;-\alpha)\frac{1}{k!}\big(e_{\lambda}(t)-1\big)^{k} \label{42}\\
	&=\ \sum_{k=0}^{\infty}C_{k}(\beta-\alpha;-\alpha)\sum_{n=k}^{\infty}S_{2,\lambda}(n,k)\frac{t^{n}}{n!} \nonumber \\
	&=\ \sum_{n=0}^{\infty}\bigg(\sum_{k=0}^{n}C_{k}(\beta-\alpha;-\alpha)S_{2,\lambda}(n,k)\bigg)\frac{t^{n}}{n!}. \nonumber
\end{align}
Therefore, by \eqref{20} and \eqref{42}, we obtain the following theorem. 
\begin{theorem}
	For $n\ge 0$, we have 
	\begin{displaymath}
	\mathrm{Bel}_{n,\lambda}(\alpha,\beta)\ =\ \sum_{k=0}^{n}C_{k}(\beta-\alpha;-\alpha)S_{2,\lambda}(n,k).
	\end{displaymath}
\end{theorem}
By replacing $t$ by $\log_{\lambda}(1+t)$ in \eqref{20}, we get 
\begin{align}\label{43}
e^{\alpha t}(1+t)^{\beta-\alpha}\ &=\ \sum_{k=0}^{\infty}\mathrm{Bel}_{k,\lambda}(\alpha,\beta)\frac{1}{k!}\big(\log_{\lambda}(1+t)\big)^{k}\\
&=\ \sum_{k=0}^{\infty}\mathrm{Bel}_{n,\lambda}(\alpha,\beta)\sum_{n=k}^{\infty}S_{1,\lambda}(n,k)\frac{t^{n}}{n!}\nonumber\\
&=\ \sum_{n=0}^{\infty}\bigg(\sum_{k=0}^{n}\mathrm{Bel}_{k,\lambda}(\alpha,\beta)S_{1,\lambda}(n,k)\bigg)\frac{t^{n}}{n!}. \nonumber
\end{align}
On the other hand, 
\begin{equation}
	e^{\alpha t}(1+t)^{\beta-\alpha}\ =\ \sum_{n=0}^{\infty}C_{n}(\beta-\alpha;-\alpha)\frac{t^{n}}{n!}. \label{44}
\end{equation}
Therefore, by \eqref{43} and \eqref{44}, we obtain the following theorem. 
\begin{theorem}
	For $n\ge 0$, we have 
	\begin{displaymath}
	C_{n}(\beta-\alpha;-\alpha)\ =\ \sum_{k=0}^{n}\mathrm{Bel}_{k,\lambda}(\alpha,\beta)S_{1,\lambda}(n,k). 
	\end{displaymath}
\end{theorem}

\section{Conclusion} 

In this paper, we introduced the two variable degenerate Bell polynomials as a generalization of the single variable degenerate Bell polynomials which were previously studied under a different name of partially degenerate Bell polynomials. Then, among other things, we showed their connections with Poisson degenerate central moments and Charlier polynomials. Their details and results are as in the following.
In Theorem 1, we expressed the two variable degenerate Bell polynomials in terms of the degenerate Bell polynomials and showed they reduce to the degenerate Bell polynomials when the two arguments are the same.
Then we defined the Poisson degenerate central moments and proved that they are equal to the two variable degenerate Bell polynomials. We found an explicit expression for the two variable degenerate Bell polynomials in Theorem 3. By letting one of the arguments of the two variable degenerate Bell polynomials equal to zero, we were led to define Stirling type numbers as the coefficients in the expansion of powers of $\alpha$. Then we expressed the Stirling type numbers in terms of the degenerate Stirling numbers of the second kind and vice versa in Theorems 4 and 5. Finally, we found certain connections between the two variable degenerate Bell polynomials and Charlier polynomials involving degenerate Stirling numbers in Theorems 8 and 9. \par
There are three immediate possible applications of our results to probability, differential equations and symmetry. As to applications to probability, we saw in this paper a connection between the two variable degenerate Bell polynomials and the Poisson degenerate central moments. For the other two possible applications of our results, we let the reader refer to [12 and the references therein]. \par
As one of our future projects, we would like to continue to study degenerate versions of certain special polynomials and numbers and their applications to physics, science and engineering as well as mathematics.


\begin{thebibliography}{9}
\bibitem{1}
Brillhart, J. \emph{Mathematical Notes: Note on the Single Variable Bell Polynomials,} Amer. Math. Monthly  \textbf{74} (1967),  no. 6, 695--696.
\bibitem{2}
Carlitz, L. \emph{Some remarks on the Bell numbers,} Fibonacci Quart. \textbf{18} (1980),  
  no. 1, 66--73.
\bibitem{3}
Carlitz, L. \emph{Single variable Bell polynomials,} Collect. Math.  \textbf{14} (1962), 13–-25. 
\bibitem{4}
Carlitz, L. \emph{Degenerate Stirling, Bernoulli and Eulerian numbers,} Utilitas Math. 15 (1979), 51–-88. 
\bibitem{5}
Howard, F. T. \emph{A special class of Bell polynomials,} Math. Comp. \textbf{35} (1980), no. 151, 977--989.
\bibitem{6}
Howard, F. T. \emph{Bell polynomials and degenerate Stirling numbers,} Rend. Sem. Mat. Univ. Padova  \textbf{61} (1979), 203--219.
\bibitem{7}
Kim, D. S.; Kim, T. \emph{A note on a new type of degenerate Bernoulli numbers,} Russ. J. Math.Phys. \textbf{27} (2020), no. 2, 227--235.
\bibitem{8}
Kim, D. S.; Kim, T. \emph{On degenerate Bell numbers and polynomials,} Rev. R. Acad.  Cienc. Exactas Fis. Nat. Ser. A Mat. RACSAM  \textbf{111} (2017), no. 2, 435--446.
\bibitem{9}
Kim, J. B. \emph{On some numbers from the Bell polynomials,} Math. Japon. \textbf{18} (1973), 1--4.
\bibitem{10}
Kim, T. \emph{A note on degenerate Stirling polynomials of the second kind,} Proc. Jangjeon Math. Soc.  \textbf{11} (2017), no. 3, 319--331.
\bibitem{11}
Kim, T.; Kim, D. S.; Dolgy, D. V. \emph{On partially degenerate Bell numbers and polynomial,} Proc. Jangjeon Math. Soc. \textbf{20} (2017), no. 3, 337--345.
\bibitem{12}
Kim, T.; Kim, D. S.; Dolgy, D. V.; Kwon J. \emph{Some identities on generalized degenerate Genocchi and Euler numbers,} Informatica 2020 (31), no. 4, 42--51.
\bibitem{13}
Kim, T.; Kim, D. S.; Jang, L.-C.; Kim, H. Y. \emph{A note on discrete degenerate random variables,} Proc. Jangjeon Math. Soc. \textbf{23} (2020), no. 1, 125--135.
\bibitem{14}
Kim, T.; Kim, D. S.; Kwon, J.; Lee, H. \emph{Degenerate polyexponential functions and type 2 degenerate poly-Bernoulli numbers and polynomials,} Adv. Difference Equ.  2020, 2020:168. 
\bibitem{15}
Kim, T.; Yao, Y.; Kim, D. S.;  Jang, G.-W. \emph{Degenerate $r$-Stirling numbers and $r$ -Bell polynomials,} Russ. J. Math. Phys.  \textbf{25} (2018), no. 1, 44--58. 
\bibitem{16}
Privault, N. \emph{Generalized Bell polynomials and the combinatorics of Possion central moments,} Elect. J. Comb. \textbf{18} (2011), paper no. 54.
\bibitem{17}
Roman, S. \emph{The umbral calculus,} Pure and Applied Mathematics, 111. Academic Press, Inc. [Harcourt Brace Jovanovich, Publishers], New York, 1984.
\bibitem{18}
Roman, S.; Rota, G.-C. \emph{The umbral calculus,} Advances in Math. \textbf{27} (1978), no. 2, 95--188.
\bibitem{19}
Ross, S.M. \emph{Introduction to Probability Models,} Academic Press, Cambridge, MA, 2007.


\end{thebibliography}
\end{document}